\documentclass[a4paper,12pt]{amsart}

\usepackage{float}
\usepackage{graphicx}

\usepackage{amsmath}

\textwidth=\paperwidth \advance\textwidth by-2in
\advance\textheight by5pt

\def\lk{\operatorname{lk}}

\def\inte{\operatorname{int}}

\theoremstyle{plain}
\newtheorem{thm}{Theorem}[section]

\newtheorem{lem}{Lemma}[section]
\theoremstyle{definition}

\theoremstyle{remark}

\begin{document}

\title[Effect of mutation] {The effect of mutation on link
concordance, 3-manifolds and the Milnor invariants}
\date{June 10, 2001 (First Edition: December 20, 2000)}
\author{Jae Choon Cha}
\email{jccha\char`\@knot.kaist.ac.kr}
\address{Department of Mathematics\\
Korea Advanced Institute of Science and Technology\\
Taejon, 305--701\\
Korea}
% Instead of \subjclass[2000]{...},
% I put the above workaround for the old versions of amsart.
\def\subjclassname{\textup{2000} Mathematics Subject Classification}
\expandafter\let\csname subjclassname@1991\endcsname=\subjclassname
\expandafter\let\csname subjclassname@2000\endcsname=\subjclassname

\begin{abstract}
We study the effect of mutation on link concordance and 3-manifolds.
We show that the set of links concordant to sublinks of homology
boundary links is not closed under positive mutation.  We show that
mutation does not preserve homology cobordism classes of 3-manifolds.
A significant consequence is that there exist 3-manifolds which have
the same quantum $\mathrm{SU}(2)$-invariants but are not homology
cobordant.  These results are obtained by investigating the effect of
mutation on the Milnor $\bar\mu$-invariants, or equivalently the
Massey products.
\end{abstract}

\maketitle

\section{Introduction and main results}
\label{sec:intr-main-results}

Let $L$ be a link in~$S^3$.  An embedded 2-sphere $S$ in $S^3$ is
called a \emph{Conway sphere} for $L$ if $S$ meets $L$ transversally
at four points.  By cutting $S^3$ along $S$ and reglueing along an
orientation preserving involution on $S$ which preserves $L\cap S$
setwisely and does not fix any points in $L\cap S$, we obtain a new
link in $S^3$ with the same number of components as~$L$.  It is called
a \emph{mutant} of~$L$.  If $L$ is oriented and the orientation of $L$
is preserved, it is called a \emph{positive mutant}.

Many invariants fail to distinguish links from their (positive)
mutants.  Skein invariants including the Alexander, Jones and HOMFLY
polynomials, are preserved by mutation, although some colored versions
of them are known to distinguish some knots from positive
mutants~\cite{Morton-Traczyk:1988-1}.  $S$-equivalence classes of
knots are preserved by positive mutation (for a proof,
see~\cite{Kirk-Livingston:1999-1}).  Consequently, all invariants
derived from Seifert matrices, including the Alexander module, torsion
invariants and signatures, are also preserved.

The problem of distinguishing knots and links from (posive) mutants
\emph{up to concordance} is even harder.  Two links $L_0$ and $L_1$ in
$S^3$ are said to be \emph{concordant} if there is a proper
submanifold $C$ in $S^3\times [0,1]$ which meets $S^3\times\{i\}$ in
$L_i$ for $i=0,1$ and is an $h$-cobordism between $L_0$ and~$L_1$.  A
link concordant to a trivial link is called a \emph{slice link}.  For
knots, primary concordance invariants are extracted from the infinite
cyclic structure of complements, but they do not detect positive
mutation.  Using the secondary invariants of Casson and
Gordon~\cite{Casson-Gordon:1978-1, Casson-Gordon:1986-1}, Kirk and
Livingston proved that positive mutation alters knot concordance
classes~\cite{Kirk-Livingston:1999-3,Kirk-Livingston:1999-1}.
Recently Cochran, Orr and Teichner found new invariants that detect
non-slice knots which are indistinguishable from slice knots by any
previously known invariants~\cite{Cochran-Orr-Teichner:1999-1}.  At
the present time, few is known about the effect of mutation on the
Cochran-Orr-Teichner invariants.

The concordance theory of links is much more complicated than that of
knots.  The crucial difficulty is that no ``good'' structure on link
complements is known.  Because of this, links that admit specific
extra structures was extensively studied.  A~link with $m$ components
is called a \emph{boundary link} if there is a homomorphism of the
link group into the free group $F$ of rank $m$ that sends an element
representing the $i$-th meridian to the $i$-th generator.
Generalizing the notion of boundary links, we call a link a
\emph{homology boundary link} if it admits a surjection of the link
group onto~$F$. Equivalence classes of boundary links and homology
boundary links under appropriate concordance that preserves the extra
structures were classified by Cappell and
Shaneson~\cite{Cappell-Shaneson:1980-1}, Ko~\cite{Ko:1987-1},
Mio~\cite{Mio:1987-1}, Duval~\cite{Duval:1986-1} and Cochran and
Orr~\cite{Cochran-Orr:1994-1}.  A more general one is the \vbox
to0mm{\vss\hbox{$\widehat F$}}-structure defined
in~\cite{Levine:1989-1, Levine:1989-2}, which \emph{sublinks of
homology boundary links} admit.

A fundamental question arising in the study of the link concordance
theory is what kinds of links are concordant to boundary links,
homology boundary links or sublinks of homology boundary links, and in
particular, whether the sets of concordance classes of boundary links,
homology boundary links and sublinks of homology boundary links are
proper subsets of the next ones.  In this regard, the relationship of
geometric operations on links and the extra structures is of special
importance.  The (strong) fusion operation was studied as a tool that
changes the extra structures~\cite{Cochran:1987-1,
Cochran-Levine:1991-1, Kaiser:1992-1}.  Especially, in the early 90's,
Cochran and Orr proved that homology boundary links which are not
concordant to boundary links can be produced by the strong fusion
operation on boundary links~\cite{Cochran-Orr:1993-1}.  It is known
that a similar effect can be achieved by mutation; the author and Ko
proved that positive mutation transforms some slice links into
sublinks of homology boundary links which are not concordant to
boundary links~\cite{Cha-Ko:2000-1,Cha-Ko:2000-2}.  One of our main
results is that the effect of mutation on link concordance classes is
much more drastic.

\begin{thm}
\label{thm:main}
There exist slice links which are transformed by positive mutation
into links not concordant to sublinks of homology boundary links.
\end{thm}

As a consequence, the set of links concordant to sublinks of homology
boundary links is not closed under positive mutation.

Mutation on knots and links is closely related with the following
analogue for 3-manifolds.  By cutting a 3-manifold $M$ along an
embedded surface $V$ of genus two and reglueing along an orientation
preserving involution on $V$ with orbit space~$S^2$, we obtain a new
manifold, which is called a \emph{mutant} of~$M$.  It is well known
that the complement and the surgery manifold of a mutant of a
\emph{knot} $K$ are obtained by performing mutation on those of $K$
(at most twice).  For the analogue for links, we introduce additional
terminologies.  Mutation on a link $L$ along a Conway sphere $S$ is
called \emph{self-mutation} if $S$ meets exactly one component of~$L$.
Otherwise, it is called \emph{bi-mutation}.  Then the above statement
about the complement and the surgery manifold is true for
self-mutation of a link, and bi-mutation of a link preserving each
component (called type $R$ in Section~\ref{sec:bi-mutation}).

It is also hard to detect the effect of mutation on 3-manifolds.  Many
invariants of 3-manifolds are preserved by mutation.  They include the
hyperbolic volume~\cite{Ruberman:1987-1}, the Reidemister
torsion~\cite{Porti:1997-1}, the $\eta$-invariant (and the
Chern-Simons invariant)~\cite{Meyerhoff-Ruberman:1990-1}, the Floer
homology~\cite{Ruberman-1999:1,Kawauchi:1996-2}, the Casson-Walker
invariant~\cite{Kirk:1989-1} and several quantum
invariants~\cite{Kania-Bartoszynska:1993-1, Kawauchi:1994-1,
Lickorish:1993-1, Rong:1994-1}.  Examples of mutation altering the
diffeomorphism type of 3-manifolds are found by
Kirk~\cite{Kirk:1989-1}, Lickorish~\cite{Lickorish:1993-1},
Kania-Bartoszynska~\cite{Kania-Bartoszynska:1993-1} and
Kawauchi~\cite{Kawauchi:1994-1}.  

We propose a new technique to detect the effect of mutation on
3-manifolds, using the lower central series.  For a group~$G$, the
lower central series $G_q$ is defined inductively by $G_1=G$,
$G_{q+1}=[G,G_q]$.  We call the quotient groups $G/G_q$ \emph{lower
central quotients}.

\begin{thm}
\label{thm:mutation-not-preserving-lcq}
The lower central quotients of the fundamental groups of 3-manifolds
are not preserved by mutation.
\end{thm}

It is well known that the lower central quotients are invariant under
\emph{homology cobordism}.  From this we obtain significant
consequences of Theorem~\ref{thm:mutation-not-preserving-lcq};
homology cobordism classes of 3-manifolds are not preserved by
mutation, and hence there are 3-manifolds which are not homology
cobordant but have the same quantum $\mathrm{SU}(2)$ invariants.

The above results on link concordance and 3-manifolds are proved by
investigating the effect of mutation on the Milnor
$\bar\mu$-invariants, or equivalently the Massey products.  Throughout
this paper, we use Milnor's notation for the invariants.  For an
oriented ordered link $L$ with $m$ components and a sequence
$I=i_1\cdots i_q$ which consists of integers between $1$ and~$m$
(denoting components of~$L$), Milnor defined an integer $\mu_L(I)$
which depends on a choice of meridian elements in~$\pi_1(S^3-L)$, and
proved that the residue class $\bar\mu_L(I)$ of $\mu_L(I)$ modulo an
integer $\Delta_L(I)$ is an isotopy invariant
of~$L$~\cite{Milnor:1957-1}.  (See Section~\ref{sec:bi-mutation} for
details.)  We denote $q$ by $|I|$ and call it the \emph{weight}
of~$\bar\mu_L(I)$.  Stallings proved that $\bar\mu_L(I)$ is an
invariant under $I$-equivalence, and hence, under link
concordance~\cite{Stallings:1965-1}.  Turaev~\cite{Turaev:1976-1},
Porter~\cite{Porter:1980-1} and Stein~\cite{Stein:1990-1} proved that
the $\bar\mu$-invariants are equivalent to the Massey cohomology
products associated to specific defining systems in link complements.

In Section~\ref{sec:bi-mutation}, the effect of bi-mutation on
$\bar\mu$-invariants is fully understood as follows.  First, we show
that bi-mutation along a Conway sphere $S$ preserves $\bar\mu(I)$ (up
to indeterminancy) if any component disjoint to $S$ is involved
in~$I$.

\begin{thm}
\label{thm:bi-mutation-preserving-mu-inv}
Suppose that $L^*$ is obtained from a link $L$ by bi-mutation along a
Conway sphere which meets the $i$, $j$-th components of $L$, and
$I=i_1\cdots i_q$ is a sequence such that $i_k \not\in \{i,j\}$ for
some $k$. Then $\mu_L(I) \equiv \mu_{L^*}(I) \mod
\gcd(\Delta_L(I),\Delta_{L^*}(I))$ for some orientation and order of
$L^*$.
\end{thm}

We remark that the orientation and the order for the mutant are left
ambiguous in Theorem~\ref{thm:bi-mutation-preserving-mu-inv} because
there is no natural way to choose them in general.  This is not a real
problem, since it causes to the $\bar\mu$-invariants the mere
ambiguity of sign and numbering of components.  In particular, it is
irrelevant to the vanishing of the $\bar\mu$-invariants.  Indeed in
the proof of Theorem~\ref{thm:bi-mutation-preserving-mu-inv} we will
clearly indicate how to choose an orientation and an order satisfying
the conclusion.

If $I$ does not satisfy the hypothesis of
Theorem~\ref{thm:bi-mutation-preserving-mu-inv}, we may forget
components disjoint to the Conway sphere and may assume that $L$ is a
2-component link without any loss of generality.  In this case we show
that higher weight $\bar\mu$-invariants are not preserved by
bi-mutation.

\begin{thm}
\label{thm:bi-mutation-not-preserving-mu-inv}
\begin{enumerate}
\item
All $\bar\mu$-invariants of weight $<6$ for $2$-component links are
preserved (up to sign) by bi-mutation.
\item
For any even $q\ge 6$, there exists a $2$-component slice link $L$
(whose $\bar\mu$-invariants vanish automatically) with a positive
bi-mutant that has vanishing $\bar\mu$-invariants of weight~$<q$ but
nonvanishing $\bar\mu$-invariants of weight~$q$.
\end{enumerate}
\end{thm}

As a consequence of the proof of
Theorem~\ref{thm:bi-mutation-not-preserving-mu-inv}, it is also shown
that positive mutation does not preserves the Cochran
invariants~\cite{Cochran:1985-2}, which are known to be integer-valued
liftings of the specific $\bar\mu$-invariants of the form
$\bar\mu(112\cdots2)$.

In Section~\ref{sec:self-mutation}, we show similar results for
self-mutation.

\begin{thm}
\label{thm:self-mutation-not-preserving-mu-inv}
\begin{enumerate}
\item
All $\bar\mu$-invariants of weight $<6$ for 2-component links are
preserved (up to sign) by self-mutation.
\item
There is a ribbon link $L$ with a positive self-mutant which has
nontrivial $\bar\mu$-invariants.
\end{enumerate}
\end{thm}

Theorem~\ref{thm:main} follows the second conclusions of
Theorems~\ref{thm:bi-mutation-not-preserving-mu-inv}
and~\ref{thm:self-mutation-not-preserving-mu-inv}, since sublinks of
homology boundary links have vanishing $\bar\mu$-invariants. 

In Section~\ref{sec:mutation-3-mfd}, we prove
Theorem~\ref{thm:mutation-not-preserving-lcq} by investigating surgery
manifolds of links.  We construct mutative pairs of 3-manifolds with
the different ninth lower central quotients.

\subsection*{Acknowledgements.}
The author would like to thank Ki Hyoung Ko and Won Taek Song for
helpful conversations.

\section{Bi-mutation and the $\bar\mu$-invariants}
\label{sec:bi-mutation}

We begin with the definition of the $\bar\mu$-invariants due to
Milnor~\cite{Milnor:1954-1,Milnor:1957-1}.  Let $L$ be an oriented
ordered link with $m$-components.  We denote the link group
$\pi_1(S^3-L)$ by~$G_L$.  Milnor proved that for any homomorphism of
the free group $F$ on $m$ generators $x_1,\ldots,x_m$ into $G_L$ which
sends $x_i$ to an element representing the $i$-th meridian, the
induced map $F/F_q \to G_L/(G_L)_q$ is surjective for all~$q$, and
furthermore, its kernel is the normal subgroup generated by
$[x_1,w_1],\ldots,[x_m,w_m]$, where $w_i$ is an arbitrary word sent to
an element representing the $i$-th longitude in~$G_L/(G_L)_q$.
Consider the Magnus expansion of $F$ into the ring of formal integral
power series in noncommutative $m$ variables $X_1, \ldots, X_m$, which
is defined by $x_i \to 1+X_i$ and $x_i^{-1} \to 1-X_i+X_i^2-\cdots$.
The Magnus expansion of $w_j$ is of the form $1+\sum \mu_L(i_1\cdots
i_k j) X_{i_1} \cdots X_{i_k}$, where the sum runs over all $i_1\cdots
i_k$.  Let $\Delta_L(i_1\cdots i_k j)$ be the greatest common divisor
of $\mu_L(J)$ where $J$ runs over all sequences obtained by permuting
proper subsequences of $i_1\cdots i_k j$ cyclically.  Milnor proved
that if $q>k+1$, the residue class $\bar\mu_L(i_1\cdots i_k j)$ of
$\mu_L(i_1\cdots i_k j)$ modulo $\Delta_L(i_1\cdots i_k j)$ is an
isotopy invariant of~$L$.  We remark that $\bar\mu(I)$ is defined for
any sequence~$I$ since an arbitrarily large $q$ can be chosen.

We will prove Theorem~\ref{thm:bi-mutation-preserving-mu-inv} using
this algebraic definition of $\bar\mu$-invariants. Suppose $S$ is a
Conway sphere of $L$ which meets two components of~$L$.  Reordering
the components, we may assume that $S$ meets the first two components
of~$L$.  Suppose that $I=i_1\cdots i_q$, where $i_k\ne 1, 2$ for some
$k$.  By the cyclic symmetry of
$\bar\mu$-invariants~\cite{Milnor:1957-1} and by reordering, we may
assume that $i_q = 3$.  $S$~bounds two 3-balls in $S^3$.  Let $B_1$ be
the one that contains the third component of~$L$, and $B_2$ be the
other.  By reordering again, we may assume that $B_2$ contains the
last $(m-n)$ components and $B_1$ contains the remaining components
except the first two components.

Let $F$ be the free group of rank $m$ as before, and let $E$ be the
subgroup of $F$ generated by the first $n$ generators
$x_1,\ldots,x_n$.  We assert that any homomorphism of $E$ into $H =
\pi_1(B_1-L)$ that sends $x_i$ to an element representing the $i$-th
meridian induces a surjection of $E/E_q$ onto~$H/H_q$.  It can be
proved by an argument similar to Milnor's proof of the analogous
result for link groups~\cite{Milnor:1957-1}, however, we give a
simpler proof using Stallings' theorem~\cite{Stallings:1965-1} as
follows.  For $i=3,\ldots,n$, choose an arc $\gamma_i$ in $B_1$
joining a point on the $i$-th component of $L$ and a point on~$S$.  We
may assume that $\gamma_i$'s are pairwisely disjoint and the interior
of each $\gamma_i$ is disjoint to~$L$.  Let $X= B_1-(L\cup
\gamma_3\cup \cdots \cup \gamma_m)$.  Then $H$ is a quotient group
of~$\pi_1(X)$, and hence the given meridian homomorphism of $E$ into
$H$ is lifted to a homomorphism of $E$ into~$\pi_1(X)$.  Since
$H_1(X)$ is a free abelian group generated by meridians and
$H_2(X)=0$, an isomorphism of $E/E_q$ onto $\pi_1(X)/\pi_1(X)_q$ is
induced by Stallings' Theorem.  The assertion follows.

By the assertion, there is a word $w$ in the first $n$ generators
$x_1,\ldots,x_n$ (and their inverses) which represents the third
longitude in $H/H_q$.  Since the composition $E\to H \to G_L$ is
extended to a homomorphism of $F$ into $G_L$ sending generators to
meridians, $w$ can also be viewed as a word representing the third
longitude in~$G_L/(G_L)_q$.  On the other hand, the same argument
works for a mutant $L^*$ of $L$ which is obtained by a mutation along
$S$, provided an orientation and an order of $L^*$ are chosen so that
the restrictions of them on $L^*\cap B_1$ coincide with those of $L
\cap B_1$.  Thus the same word $w$ represents the third longitude of
$L^*$ in $G_{L^*}/(G_{L^*})_q$.  This completes the proof of
Theorem~\ref{thm:bi-mutation-preserving-mu-inv}.

\begin{figure}[hbt]
%\fboxsep=10mm \fbox{\texttt{Connected sum decomposition of $L$}}
\includegraphics{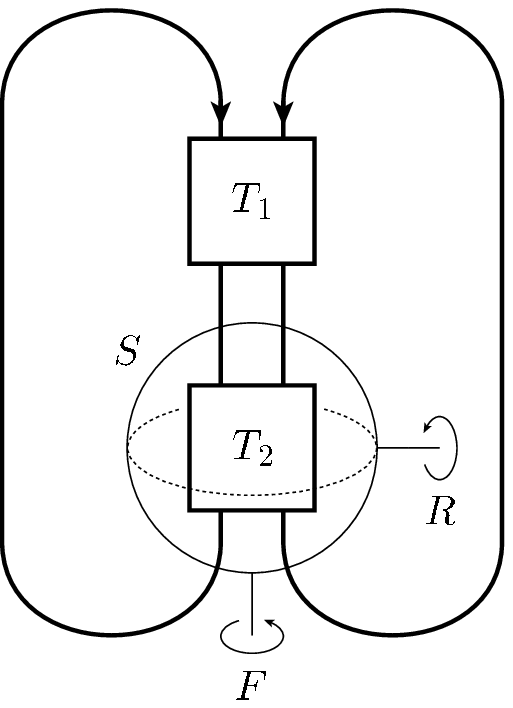}
\caption{}
\label{fig:csum-decomp}
\end{figure}

The rest of this section is devoted to the study of the effect of
bi-mutation on the $\bar\mu$-invariants of two component links.
Suppose that $L$ is an oriented ordered two component link and $S$ is
a Conway sphere which intersects both components of~$L$ as shown in
Figure~\ref{fig:csum-decomp}.  An involution on $S$ that produces a
mutant of $L$ is isotopic to either one of the followings: the
$\pi$-rotations along the axes $F$ and $R$ shown in
Figure~\ref{fig:csum-decomp}, or the composition of them.  The
associated bi-mutations will be called \emph{bi-mutations of type $F$,
$R$ and~$FR$}, and the associated bi-mutants will be denoted by $L^F$,
$L^R$ and~$L^{FR}$, respectively.  They can be characterized as
follows: bi-mutation of type $F$ is the positive mutation, which
preserves the orientation and reverses the order of the components
of~$L$.  Bi-mutation of type $R$ reverses the orientation and
preserves the order.  Bi-mutation of type $FR$ reverses both of the
orientation and the order.

Cutting $S^3$ along $S$, we obtain two string links denoted by $T_1$
and $T_2$ in Figure~\ref{fig:csum-decomp}.  $L$~can be viewed as a
connected sum of their closures $\alpha$ and~$\beta$.  For
convenience, we always choose an orientation and an order of a
bi-mutant which coincide with those of~$T_1$.  Note that this is
irrelevant to the vanishing of $\bar\mu$-invariants.

We will describe the $\bar\mu$-invariants of the bi-mutants in terms
of those of $\alpha$ and~$\beta$.  For this, we introduce some
operations on indices.  For a finite sequence $I$ consisting of $1$
and~$2$, let $I^F$ (resp.\ $I^R$) be the sequence obtained from $I$ by
exchanging $1$ and $2$ (resp. by reversing indices), and let $I^{FR} =
(I^F)^R$.  For example, for $I=112222$, we have $I^F = 221111$, $I^R =
222211$ and $I^{FR} = 111122$.

\begin{lem}
\label{lem:mu-inv-of-bi-mutant}
\begin{enumerate}
\item
$D(I) = \gcd(\Delta_{\alpha}(I), \Delta_{\beta}(I))$ divides
$\Delta_L(I)$ and
$$
\mu_L(I) \equiv \mu_{\alpha}(I) + \mu_{\beta}(I) \mod D(I).
$$
\item
For $\tau=F, R$ and $FR$, $D^\tau(I) = \gcd(\Delta_{\alpha}(I),
\Delta_{\beta}(I^\tau))$ divides $\Delta_{L^\tau}(I)$ and
$$
\mu_{L^\tau}(I) \equiv \mu_{\alpha}(I) + \mu_{\beta}(I^\tau) \mod D^\tau(I).
$$
\end{enumerate}
\end{lem}

\begin{proof}
Since $L$ is a connected sum of $\alpha$ and $\beta$, the first
conclusion follows the additivity of the $\bar\mu$-invariants under
connected sum (e.g. see~\cite{Krushkal:1998-1}).

Fix $q>|I|$, and let $w_i(x_1,x_2)$ be a word representing the $i$-th
longitude of $\beta$ in $G_\beta/(G_\beta)_q$.  For $L^F$, the link
$\beta^F$ obtained by reversing the order of $\beta$ plays the role
of~$\beta$.  $w_{3-i}(x_2,x_1)$ represents the $i$-th longitude
of~$\beta^F$.  Thus $\bar\mu_{\beta^F}(I) = \bar\mu_\beta(I^F)$, and
the conclusion for $L^F$ follows the additivity again.

For $L^R$, the link $\beta^R$ obtained by reversing the orientation of
$\beta$ plays the role of~$\beta$.  In this case, the word
$v_i(x_1,x_2)$ obtained by ``reversing'' the word $w_i(x_1,x_2)$ (that
is, reading from the end to the beginning) represents the $i$-th
longitude of~$\beta^R$.  Thus for $I=i_1\cdots i_k j$,
$\bar\mu_{\beta^R}(I) = \bar\mu_{\beta}(i_k\cdots i_1 j) =
\bar\mu_\beta(I^R)$ by the cyclic symmetry.  This proves the
conclusion for~$L^R$.

The conclusion for $L^{FR}$ follows the above results, since the
mutation of type $FR$ is the composition of the mutations of type $F$
and~$R$.
\end{proof}

Now we are ready to prove the invariance of lower weight
$\bar\mu$-invariants under bi-mutation.  For $2$-component links, it
is well known that $\bar\mu(12)$ and~$\bar\mu(1122)$ are all of the
nontrivial $\bar\mu$-invariants of weight~$<6$.  $\bar\mu(12)$ is the
linking number, which is preserved by mutation (up to sign because of
the ambiguity of an orientation).  The indeterminancy $\Delta(1122)$
is equal to the linking number, and hence preserved.  We may assume
that $\alpha$ has the same linking number as $L$ and $\beta$ has
linking number zero, by putting additional full twists on $\alpha$ and
$\beta$ which are cancelled in $L$ and~$L^\tau$.  Then we have
$\Delta_\alpha(1122)=\lk(L)$ and $\Delta_\beta(1122)=0$ so that
$D(1122) = \Delta_L(1122) = \lk(L) = \Delta_{L^\tau}(1122) =
D^\tau(1122)$.  By the cyclic symmetry of $\bar\mu$-invariants,
$\bar\mu_\beta(1122)=\bar\mu_\beta(1122^\tau)$ for any~$\tau$.  Thus
$\bar\mu_L(1122) = \bar\mu_{L^\tau}(1122)$ by
Lemma~\ref{lem:mu-inv-of-bi-mutant}.  This proves the first conclusion
of Theorem~\ref{thm:bi-mutation-not-preserving-mu-inv}.

Next, we will show that higher weight $\bar\mu$-invariants are not
preserved by bi-mutation.  Explicitly, for any given type $\tau$ we
will construct a link $L$ with the following properties, by choosing
$T_1$ and $T_2$ carefully.  (1)~$L$~is a ribbon link (in particular,
all $\bar\mu$-invariants of $L$ vanish).  (2)~For some $I$, all
$\bar\mu$-invariants of $L^\tau$ of weight~$<|I|$ vanish but
$\bar\mu_{L^\tau}(I)$ does not vanish.

Take the mirror image of~$T_1$ (with respect to a horizontal mirror)
as~$T_2$.  Then $\beta$ is the mirror image $\alpha^{-1}$ of $\alpha$.
Since $L$ is a connected sum of $\alpha$ and its mirror image~$\beta$,
$L$ is a ribbon link.  For any $\tau=F, R$ and $FR$,
\begin{equation*}
\begin{split}
\bar\mu_{L^\tau}(I) &\equiv
\bar\mu_{\alpha}(I)+\bar\mu_{\beta}(I^\tau) \\
&\equiv
\mu_{\alpha}(I)-\bar\mu_{\alpha}(I^\tau) \mod D(I)
\end{split}
\end{equation*}
by Lemma~\ref{lem:mu-inv-of-bi-mutant} and by the fact
$\bar\mu_{\alpha^{-1}}(I) = -\bar\mu_\alpha(I)$.  Thus if $\alpha$ has
vanishing $\bar\mu$-invariants of weight $<q$ for some~$q$, so
does~$L^\tau$.  Furthermore, if $\bar\mu_\alpha(I) \neq
\bar\mu_\alpha(I^\tau)$ for some $I$ with $|I|=q$, then
$\bar\mu_{L^\tau}(I)$ does not vanish.  Thus we are naturally led to
the question whether the invariants $\bar\mu(I)$ and $\bar\mu(I^\tau)$
are the same (for links with vanishing $\bar\mu$-invariants of
weight~$<q$) for all $I$ with $|I|=q$.

\begin{figure}[hbt]
%\fboxsep=10mm \fbox{\texttt{Milnor's link}}
\includegraphics{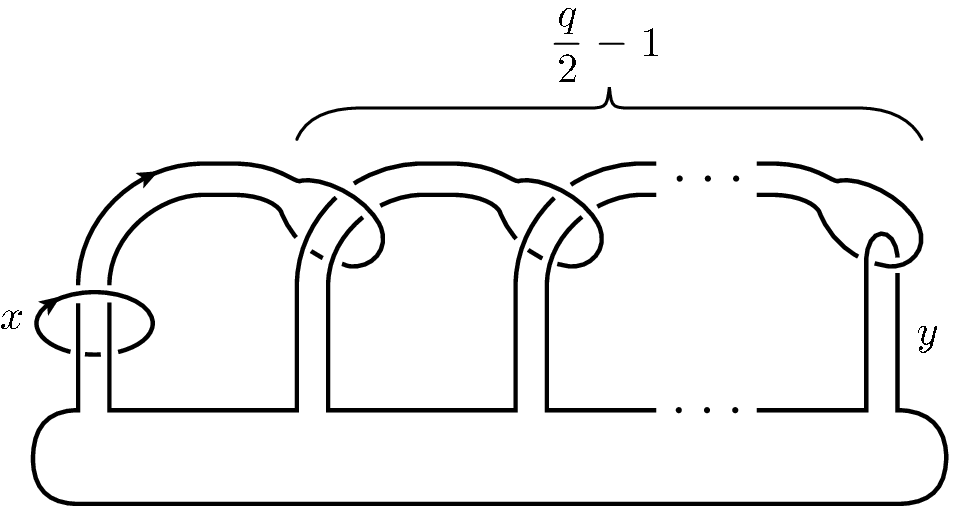}
\caption{}
\label{fig:milnor-link}
\end{figure}

In general, the answer is \emph{no} so that our construction is
successful.  A concrete example for type $F$ is obtained by taking as
$\alpha$ the link $L_q$ suggested by
Milnor~\cite[p.~301]{Milnor:1957-1} for any even $q \ge 6$.  We
illustrate $L_q$ in Figure~\ref{fig:milnor-link}, where the first and
second components are marked as $x$ and $y$, respectively.  We take as
$T_1$ any 2-string link whose closure is $\alpha$.  $L_q$~has
vanishing $\bar\mu$-invariants of weight~$<q$, and for $I=1122\cdots
2$ ($|I|=q$), $\bar\mu_{L_q}(I)=(-1)^{q/2}$ and $\bar\mu_{L_q}(I^F) =
0$~\cite{Milnor:1957-1, Cochran:1990-1}.  Thus the links $L$ and $L^F$
have the desired properties.  In particular, the second conclusion of
Theorem~\ref{thm:bi-mutation-not-preserving-mu-inv} follows.
Furthermore, $L$~and $L^{FR}$ also have the desired properties, since
$\bar\mu(I^{FR}) = \bar\mu(I^F)$ for $I=1122\cdots 2$ by the cyclic
symmetry.

For type $R$, we need to consider a more complicated example.  Take
the link given in~\cite[Figure 2.13a]{Cochran:1990-1} as $\alpha$.
Cochran proved that it has vanishing $\bar\mu$-invariants of
weight~$<10$ and nonvanishing $\bar\mu(I)$ for $I=2222121211$
in~\cite[Example 2.12]{Cochran:1990-1}.  The same computational
technique proves that it has vanishing $\bar\mu(I^R)$.  Thus a link
with the desired property for type $R$ can be produced by our
construction.  The author does not know whether there is a simpler
$\bar\mu$-invariant which is not preserved by bi-mutation of type~$R$.

We finish this section with some remarks on the problem of deciding
whether $\bar\mu(I)$ is invariant under bi-mutation of type $\tau$, or
equivalently whether $\bar\mu(I)\ne \bar\mu(I^\tau)$, for a
given~$I$. It is closely related to the (not well understood)
structure of the space $V_q$ of $\bar\mu$-invariants of
weight~$q=|I|$.  Milnor found some relations on the
$\bar\mu$-invariants~\cite{Milnor:1957-1}, however, it is not known
whether they generate all relations.  By Orr, $\dim(V_q)$ was
explicitly computed~\cite{Orr:1989-1}.  For $q<6$, $\dim(V_q)=0$ or
$1$ and actually all $\bar\mu$-invariants of weight~$<6$ are preserved
by bi-mutation, as it was shown in the above.  For most large values
of $q$, however, the space $V_q$ seems large enough to find an $I$ of
weight $q$ such that $\bar\mu(I)$ and $\bar\mu(I^\tau)$ are
independent.  On the other hand, Cochran found an algorithm which
decides whether $\bar\mu(I)$ is nontrivial for a given~$I$, (if the
answer is yes) computes the least positive value $m(I)$ of
$\bar\mu(I)$, and constructs a link
realizing~$m(I)$~\cite{Cochran:1990-1}.  Since $m(I)\ne m(I^\tau)$
implies $\bar\mu(I)\ne \bar\mu(I^\tau)$, Cochran's algorithm can be
used to obtain partial information on the decision problem.

\section{Self-mutation and the $\bar\mu$-invariants}
\label{sec:self-mutation}

In this section we investigate the effect of self-mutation on
$\bar\mu$-invariants for 2-component links.  The main tool of this
section is the dual interpretation of the Massey cohomology products
through the intersection theory of chains.  In particular, the
following geometric interpretation of Sato and Levine will be used for
studying~$\bar\mu(1122)$.  For a $2$-component link $L$ with vanishing
linking number, each component bounds a Seifert surface disjoint to
the other component.  We may assume that the two surfaces intersect
transversally along a 1-manifold.  In~\cite{Sato:1984-1}, it was shown
that the self-linking $\beta(L)$ of the 1-manifold along the framing
induced by the surfaces is a well-defined link concordance invariant.
$\beta(L)$~is called the Sato-Levine invariant.  It is well known that
$\beta(L)$ coincides with~$\bar\mu_L(1122)$.

\begin{figure}[hbt]
%\fboxsep=10mm \fbox{\texttt{Decomposition for self-mutants}}
\includegraphics{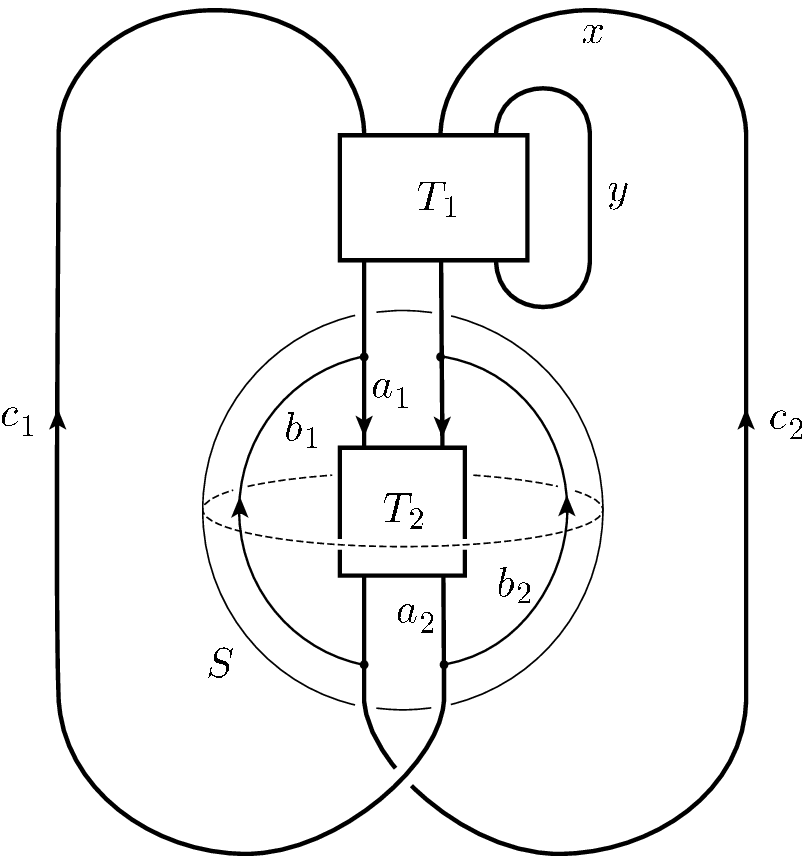}
\caption{}
\label{fig:self-mutant-decomp}
\end{figure}

Suppose that $L$ is a link with 2-components, and $L^*$ is a mutant
produced by mutation along a Conway sphere $S$ which meets only one
component of~$L$.  We will show $\bar\mu_L(1122) =
\bar\mu_{L^*}(1122)$.  First we consider a special case of a link $L$
with vanishing linking number.  For convenience, we fix an orientation
of~$L$. (Note that $\bar\mu(1122)$ is independent of the choice of an
orientation.)  Let $x$ be the component of $L$ intersecting $S$, and
$y$ be the other component.  Let $B_1$ be the 3-ball in $S^3$ which is
bounded by $S$ and disjoint to~$y$, and let $B_2=S^3-\inte B_1$.  Then
$x \cap B_1$ consists of two oriented arcs $a_1$ and~$a_2$.  Choose
disjoint oriented arcs $b_1$, $b_2$ on $S$ such that $\partial b_i =
-\partial a_i$ and $b_1\cup b_2$ is setwisely preserved by the
involution on $S$ producing~$L^*$.  $u_i=a_i\cup b_i$ is an oriented
simple closed curve.  $x \cap B_2$ consists of two oriented arcs.  Let
$c_1$ be the one joining the endpoints of $a_2$ and $a_1$, and $c_2$
be the other.  See Figure~\ref{fig:self-mutant-decomp}, where $T_1$
and $T_2$ represent a 3-string link and a 2-string link, respectively.

\def\fakeprime{{\vphantom{\prime}}}

\begin{figure}[hbt]
%\fboxsep=10mm \fbox{\texttt{The intersection of two Seifert surfaces}}
\includegraphics{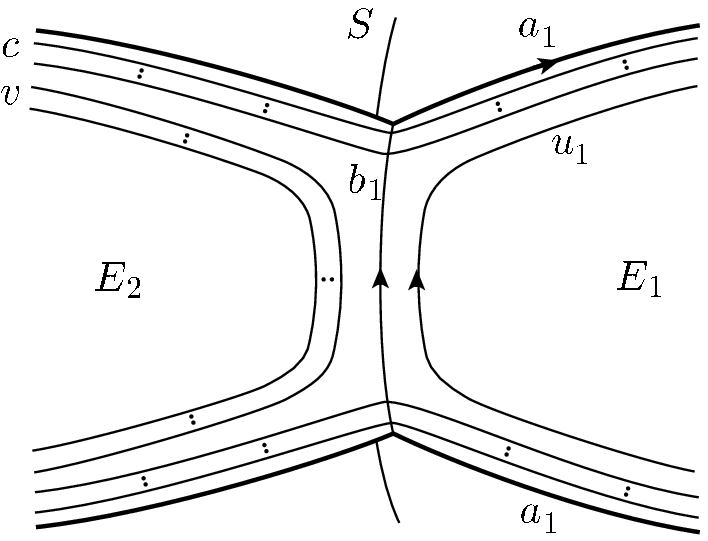}
\caption{}
\label{fig:seifert-surface-intersection}
\end{figure}

There exist oriented surfaces $E_1$ and $E_2$ properly embedded in
$B_1$ and $B_2$ which are bounded by $u_1 \cup u_2$ and $(-b_1)\cup
c_1\cup (-b_2) \cup c_2$, respectively.  Since the linking number of
$x$ and $y$ is zero, we may assume $E_2$ is disjoint to~$y$.  Then $V
= E_1\cup E_2$ is a Seifert surface of $x$ which is disjoint to~$y$.
Choose a Seifert surface $W'$ of $y$ in~$B_2$ which meets $x$
transversally.  By deleting an open tubular neighborhood of $W'\cap x$
from $W'$ and attaching thin cylinders contained in a tubular
neighborhood of $x-a_2$, a Seifert surface $W$ of $y$ which is
disjoint to $x$ is constructed.  The Sato-Levine invariant $\beta(L)$
is the self-linking of $c=V\cap W$.  We will compute it by evaluating
the value of the Seifert form of $V$ at~$(c,c)$.  $c = (E_1\cap W)
\cup (E_2 \cap W)$ and $E_1 \cap W$ is equal to the intersection of
$E_1$ and attached cylinders, which is a union of parallel copies of
$a_1$ in~$E_1$.  Therefore $c$ is homologous (on $V$) to $ku_1+v$ for
some integer $k$ and some $1$-cycle $v$ in~$E_2$.  See
Figure~\ref{fig:seifert-surface-intersection}.  Thus $\beta(L) =
\lk(c,c') = k^2 \lk(u_1',u_1^\fakeprime) + \lk(v',v) + k\lk(u_1',v) +
k\lk(u_1,v')$, where~$(\cdot)'$ denotes a cycle obtained by pushing
slightly along the positive normal direction of $V$.  Since $u_1$ and
$v$ are separated by~$S$, the last two terms vanish.  Hence $\beta(L)
= k^2 \lk(u_1',u_1^\fakeprime) + \lk(v',v)$.

Since $b_1\cup b_2$ is preserved by the involution on $S$
producing~$L^*$, a Seifert surface $V^*$ of a component of $L^*$ in
$S^3$ is obtained by glueing $(B_1,E_1)$ and $(B_2,E_2)$ along the
involution (we need to reverse the orientation of $E_1$ if the
mutation is not positive).  A Seifert surface $W^*$ of the other
component of $L^*$ such that $W^*\cap B_2 = W\cap B_2$ is obtained by
attaching cylinders to punctured $W'$.  As before, the self-linking of
$V^*\cap W^*$ is equal to $k^2\lk(u_j',u_j^\fakeprime) + \lk(v',v)$,
where $j=1$ if the mutation preserves $a_1$ setwisely, and $j=2$
otherwise.  If $j=1$, we have $\beta(L)=\beta(L^*)$ obviously.  If
$j=2$, we need additional arguments.  Since $u_1'$ and $E_1$ are
disjoint, $0=\lk(u_1', \partial E_1) = \lk(u_1',u_1^\fakeprime) +
\lk(u_1',u_2^\fakeprime)$.  Similarly $0=\lk(u_2', u_1^\fakeprime) +
\lk(u_2',u_2^\fakeprime)$. Since $u_1$ and $u_2$ are disjoint,
$\lk(u_1',u_1^\fakeprime) = -\lk(u_1',u_2^\fakeprime) =
-\lk(u_2',u_1\fakeprime) = \lk(u_2',u_2^\fakeprime)$.  This proves
$\beta(L)=\beta(L^*)$.

For the general case, choose a link $L'$ such that $\lk(L')=-\lk(L)$.
Choosing basings for connected sums carefully, a connected sum $L\#
L'$ can be viewed as a self-mutant of a connected sum $L^*\# L'$.
Since $L \# L'$ has vanishing linking number, $\mu_{L\#
L'}(1122)=\mu_{L^*\# L'}(1122)$.  By the additivity of
$\bar\mu$-invariants, $\mu_L(1122)+\mu_{L'}(1122) \equiv
\mu_{L^*}(1122)+\mu_{L'}(1122) \mod \lk(L)$.  Thus
$\bar\mu_L(1122)=\bar\mu_{L^*}(1122)$.  This completes the proof of
the first part of
Theorem~\ref{thm:self-mutation-not-preserving-mu-inv}.

\begin{figure}[hbt]
% \fboxsep=10mm \fbox{\texttt{A slice link}}
\includegraphics{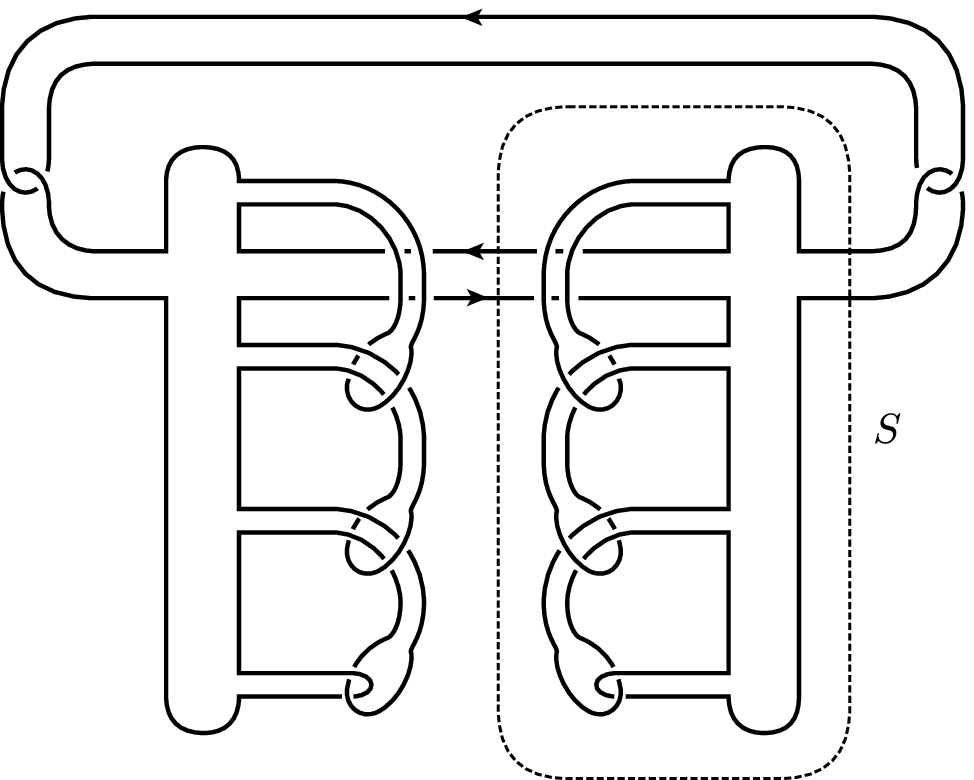}
\caption{}
\label{fig:self-mutation-ex}
\end{figure}

The rest of this section is devoted to the proof of the second part of
Theorem~\ref{thm:self-mutation-not-preserving-mu-inv}.  Let $L$ be the
link shown in Figure~\ref{fig:self-mutation-ex}.  Since $L$ is a connected
sum of a link with its mirror image, $L$ is a ribbon link.

Let $L^*$ be a mutant of $L$ obtained by the positive self-mutation
along the Conway sphere $S$ shown in
Figure~\ref{fig:self-mutation-ex}.  We will compute the
$\bar\mu$-invariants of $L^*$ using a geometric method suggested by
Cochran~\cite{Cochran:1990-1}.  In the below we describe a
reformulated version of Cochran's method where necessary surfaces are
constructively obtained.

Suppose that an oriented link with components $x_1,x_2,\ldots$ such
that $\lk(x_i,x_j)=0$ for $i\ne j$ is given.  \emph{Brackets} in
symbols $x_1,x_2,\ldots$ are defined inductively as follows.
$x_1,x_2,\ldots$ are called \emph{$1$-brackets}.  For a $k$-bracket
$u$ and an $(n-k)$-bracket $v$, $w=(u,v)$ is called an
\emph{$n$-bracket}.  $n$~is called the \emph{weight} of~$w$ and
denoted by~$|w|$.  A \emph{formal $n$-linking} is defined to be an
abstract symbol $\epsilon \lk(w)$, where $w$ is a bracket with weight
$n>1$ and $\epsilon=+$~or~$-$.  An equivalence relation on formal
linkings is generated by two relations $\pm\lk((u,v),w) \sim
\pm\lk(u,(v,w))$ and $\pm\lk(w) \sim \mp\lk(w')$ where $w'$ is
obtained by replacing a proper (i.e.\ $w\ne(u,v)$) sub-bracket $(u,v)$
in $w$ by~$(v,u)$.  If a formal $n$-linking $\pm\lk(u,v)$ has minimal
$\big||u|-|v|\big|$ in its equivalence class, it is called a
\emph{minimal $n$-linking}.

We will choose oriented surfaces $V(w)$ and oriented closed
1-manifolds $c(w)$, $c'(w)$ associated to brackets, from which the
values of formal linkings are defined.  For convenience, we will keep
the condition that if $w$ and $w'$ are as in the second relation in
the above, then the curves associated to $w$ and $w'$ are the same
with opposite orientations, and if $w'$ is obtained by replacing a
(possibly non-proper) sub-bracket $(u,v)$ in $w$ by $(v,u)$, then the
surfaces associated to $w$ and $w'$ are the same with opposite
orientations.  For $1$-brackets, let $c(x_i)=x_i$ and choose a
$0$-linking parallel $c'(x_i)$ of~$x_i$ lying on the boundary of a
tubular neighborhood $U$ of the given link.  Since
$\lk(c(x_j),c'(x_i))=0$, we can choose a Seifert surface $V(x_i)$ of
$c'(x_i)$ such that $V(x_i)\cap U=c'(x_i)$ and $V(x_i)$ and $V(x_j)$
intersect transversally for all $i,j$.  Let $c(x_i,x_j)$ be the
oriented intersection of $V(x_i)$ and $V(x_j)$ (so that the triple of
a positive tangent vector of $c(x_i,x_j)$ and the positive normal
vectors of $V(x_i)$ and $V(x_j)$ induces the orientation of the
ambient space) and choose a parallel $c'(x_i,x_j)$ of $c(x_i,x_j)$
along the framing induced by the surfaces.  We repeat this process as
follows.  Suppose that $V(u)$, $c(u,v)$ and a parallel $c'(u,v)$ of
$c(u,v)$ have been chosen for $|u|,|v|<n$ so that the followings are
satisfied.
\begin{enumerate}
\item
All surfaces are in general position.
\item
$\partial V(u)=c'(u)$ for $|u|<n$.
\item
$c(u,v)$ is the oriented intersection of $V(u)$ and $V(v)$ for
$|u|,|v|<n$.
\item
$c'(u,v)$ is a parallel on the boundary of a sufficiently small
tubular neighborhood $N(u,v)$ of $c(u,v)$ and both $V(u)$ and $V(v)$
are disjoint to~$c'(u,v)$ for $|u|,|v|<n$.
\item\label{item:closedness}
Unless $V(v)$ contains $c(u)$ or $c'(u)$, $N(u)$ is disjoint to $V(v)$
for $|u|,|v|<n$.
\end{enumerate}
In particular, $c(u)$ and $c'(u)$ have been chosen for $|u|\le n$.  As
a convention, $c(u,u)$ and $c'(u,u)$ are defined to be empty.  Note
that for $u\ne v$, $V(u)\cap V(v)$ is closed if and only if $c'(u)
\cap V(v)=V(u) \cap c'(v)=\emptyset$.  If $\lk(c(u),c'(v))=0$ for any
$|u|,|v|\le n$, we can choose $V(u)$ bounded by $c'(u)$ for $|u|=n$
(and modify $V(u)$ for $|u|<n$ if necessary) so that
(\ref{item:closedness}) is satisfied for $|u|,|v| \le n$.  Now
$c(u,v)=V(u) \cap V(v)$ is a closed 1-manifold for $|u|,|v|\le n$, and
$c'(u,v)$ can be chosen appropriately. This completes the construction
of $V(u)$, $c(u,v)$ and $c'(u,v)$ for $|u|,|v|<n+1$.

We continue this process until $\lk(c(u),c(v))$ is nonzero for some
brackets $u$ and $v$ such that $|u|\le |v|=n$.  At this time we have
surfaces $V(w)$ for $|w|\le n-1$ and curves $c(w)$, $c'(w)$ for
$w=(u,v)$, $|u|,|v| \le n-1$.  Suppose~$q\le 2n$.  Then
in~\cite[Proposition~2.11]{Cochran:1990-1} it was shown that $c(u)$
and $c(v)$ exist for any minimal $q$-linking~$\lk(u,v)$.  We define
the value of a formal $q$-linking $\lk(w)$ to be $\pm\lk(c(u),c'(v))$,
where $\pm\lk(u,v)$ is a minimal $q$-linking equivalent to~$\lk(w)$.
If all linkings of weight~$<q$ vanish, the collection of surfaces and
curves is called a \emph{surface system of weight~$q$}.  An argument
in~\cite{Cochran:1990-1} shows the following.

\begin{thm}\label{thm:minimal-linking}
Suppose that there exists a surface system of weight~$q$.  Then
\begin{enumerate}
\item
All $\bar\mu$-invariants of weight~$<q$ vanish.
\item
For $I=i_1\cdots i_{q-1} i_q$ with $i_1\ne i_q$,
$$
\bar\mu(I) = (-1)^q\sum_w \lk(w,x_{i_q})
$$
where the sum runs over all binary parenthesizations on the string
$x_{i_1}\cdots x_{i_{q-1}}$.
\end{enumerate}
\end{thm}

The sum is indeed a linear combination of minimal $q$-linkings, since
each summand is equal to a minimal $q$-linking up to sign.
Practically, the sum can be computed using the technique of the
``formal Massey products''~\cite{Cochran:1990-1}.  For example, a
straightforward computation shows that if a link with two components
$x$ and $y$ admits a surface system of weight~$9$,
$$
\bar\mu(122121222) = -20\lk(yyxy,yxyxy) -20\lk(yyxy,yyxxy)
-20\lk(yyxy,(yxy,xy))
$$
where $a_1\cdots a_n$ denotes the bracket
$(a_1,(\cdots,(a_{n-2},(a_{n-1},a_n))\cdots))$.

\begin{figure}[H]
\begin{center}
\includegraphics{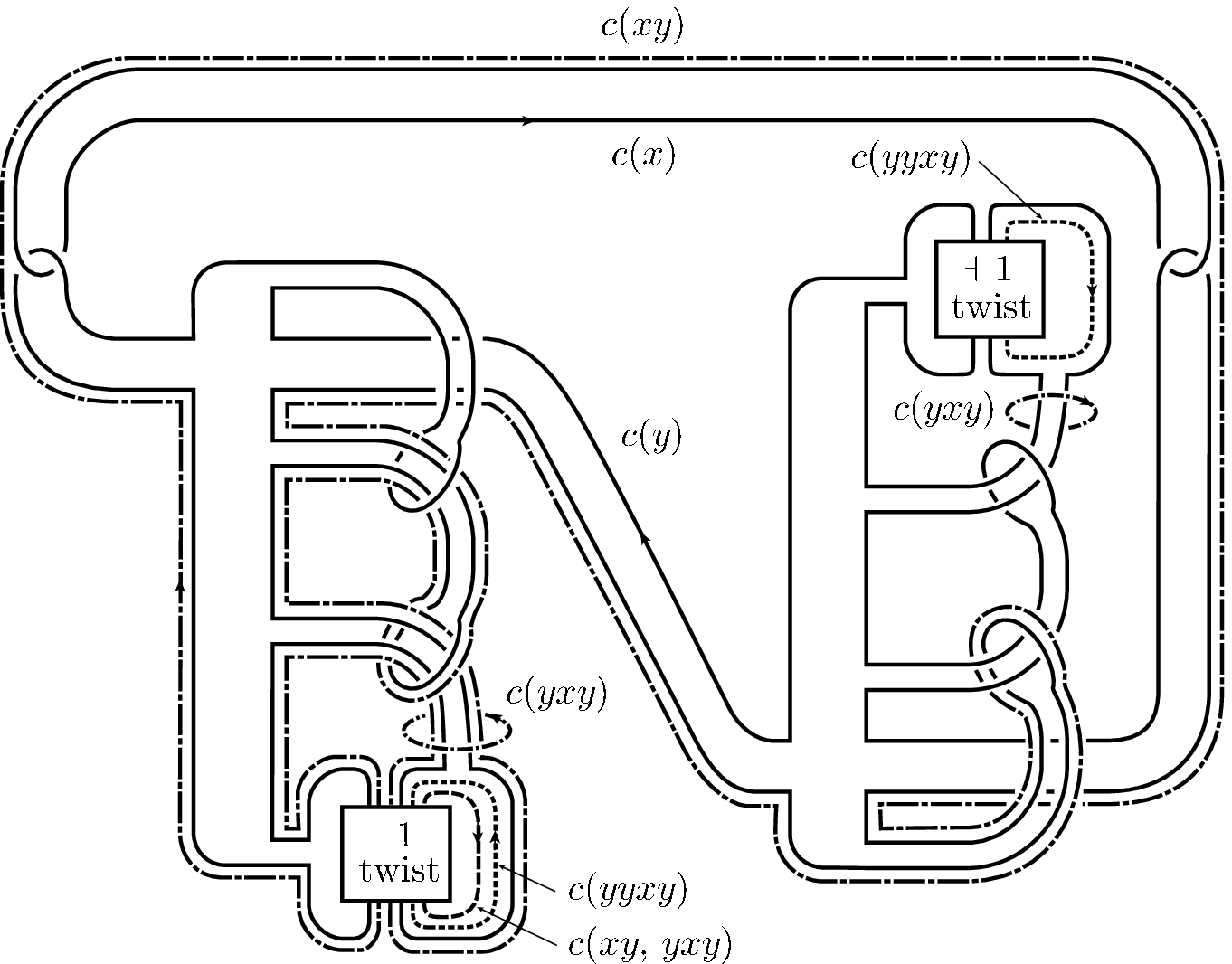}
\end{center}
\caption{}
\label{fig:curves}
\end{figure}

\begin{figure}[H]
$$V(x)=\quad\vcenter{\hbox{\includegraphics{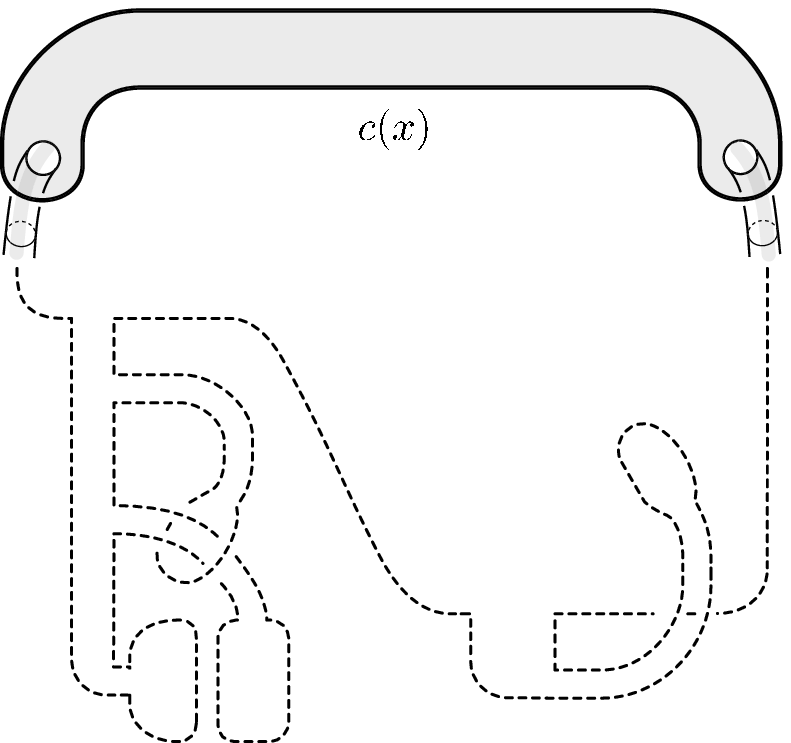}}}$$
\caption{}
\label{fig:surface-x}
\end{figure}

\begin{figure}[H]
$$V(y)=\quad\vcenter{\hbox{\includegraphics{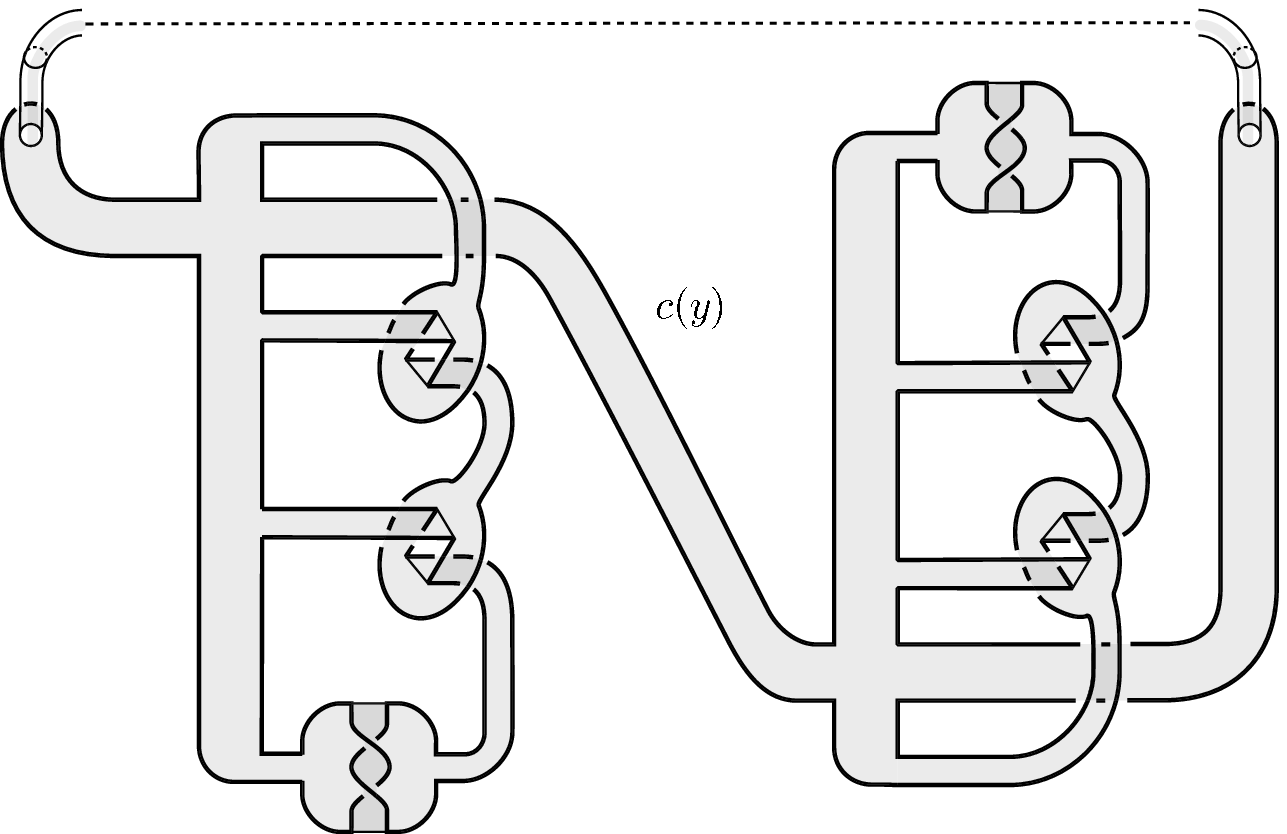}}}$$
\caption{}
\label{fig:surface-y}
\end{figure}

\begin{figure}[H]
$$V(xy)=\quad\vcenter{\hbox{\includegraphics{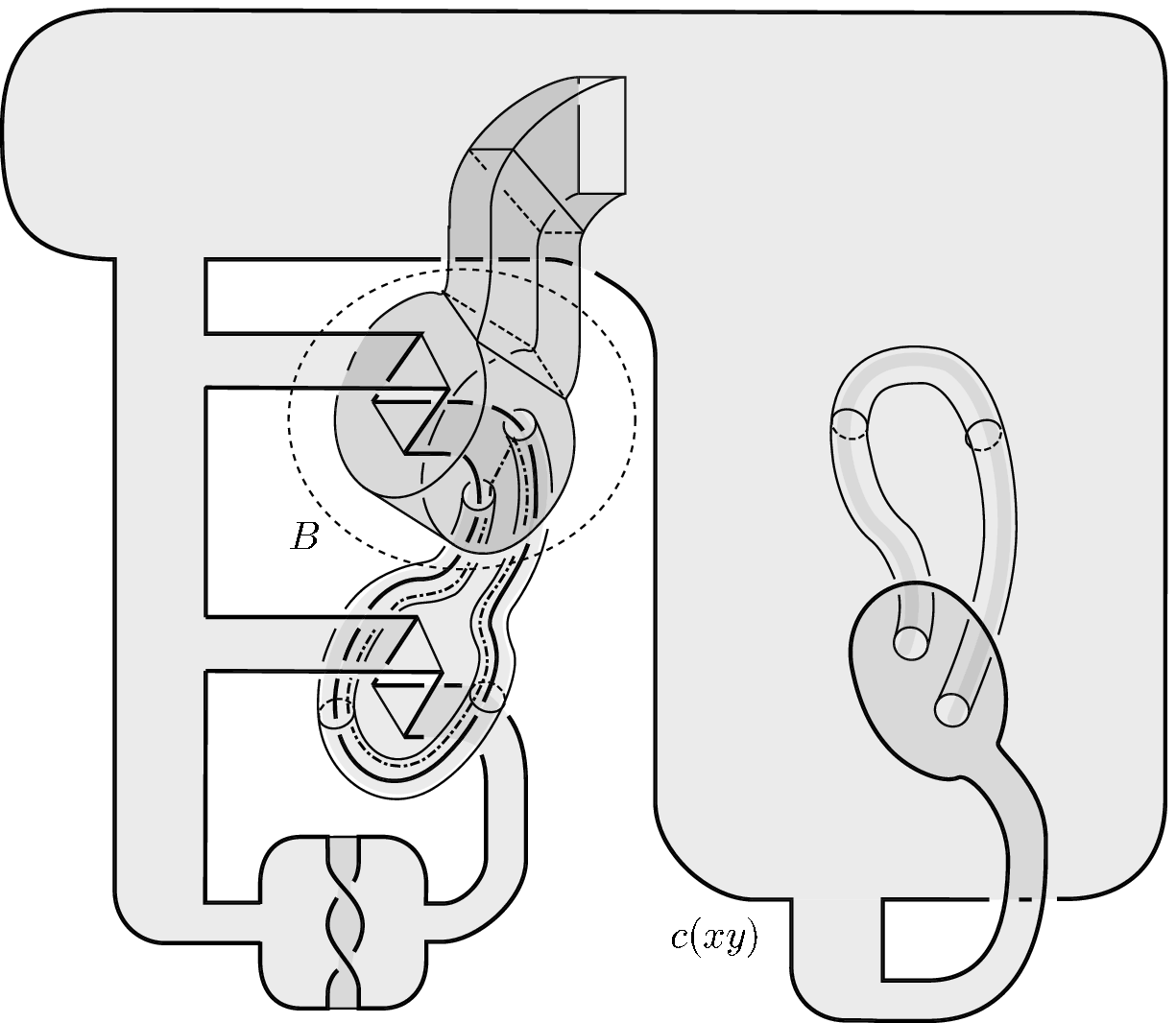}}}$$
\caption{}
\label{fig:surface-xy}
\end{figure}

\begin{figure}[H]
\begin{center}
\includegraphics{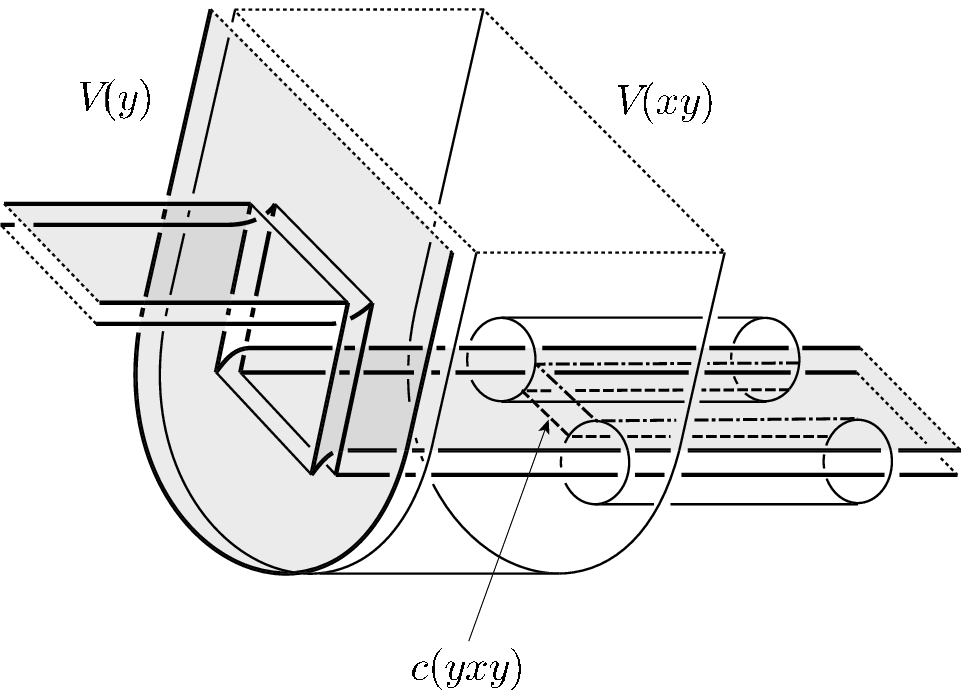}
\end{center}
\caption{}
\label{fig:surface-y-xy}
\end{figure}

\begin{figure}[H]
$$V(yxy)=\quad\vcenter{\hbox{\includegraphics{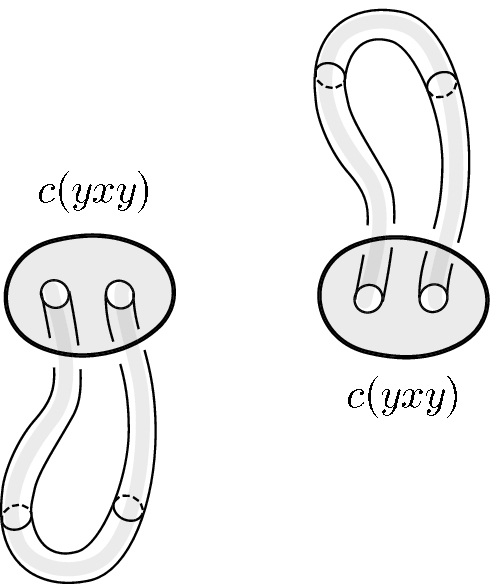}}}$$
\caption{}
\label{fig:surface-yxy}
\end{figure}

We apply this method to our example.  The link $L^*$ admits a surface
system of weight~$9$.  It is best illustrated by figures.  Nonempty
curves $c(w)$ for $|w|\le 5$ and surfaces $V(w)$ for $|w| \le 3$ are
illustrated in Figure~\ref{fig:curves} and
Figures~\ref{fig:surface-x}--\ref{fig:surface-yxy}, respectively.  For
$|w|=4$, surfaces $V(w)$ can be constructed as in the above discussion
since $\lk(c(u),c'(v))=0$ for $|u|,|v| \le 4$.  We remark that $V(xy)$
in Figure~9 is actually an immersed surface with a self-intersection
marked by~$\vcenter{\hbox{\vrule width 4pt height .8pt \hskip 1pt
\vrule width 1pt height .8pt
\hskip 1pt \vrule width 4pt height .8pt \hskip 1pt \vrule width 1pt
height .8pt \hskip 1pt \vrule width 4pt height .8pt}}$.  This is
irrelevant in computing the Massey products.  Indeed, to avoid the
self-intersection, one may take as $V(xy)$ the embedded orientable
surface obtained by splicing the immersed surface along the
self-intersection.  The most complicated part of our configuration of
surfaces is the relative location of $V(y)$ and~$V(xy)$.  For the
reader's convenience, in Figure~\ref{fig:surface-y-xy}, we illustrate
in detail parts of $V(y)$ and $V(xy)$ in the $3$-ball $B$ shown in
Figure~\ref{fig:surface-xy}.

It is straightforward to check that all minimal linkings of
weight~$\le 9$ vanish, except $\lk(yyxy, (yxy,xy))=1$.  Therefore the
$\bar\mu$-invariants of weight~$\le 8$ vanish for~$L^*$, and
$\bar\mu_{L^*}(122121222)=-20$.  This proves
Theorem~\ref{thm:self-mutation-not-preserving-mu-inv}.

We remark that for $n\ge 1$, we can construct a ribbon link and its
self-mutant with vanishing minimal linkings of weight~$\le 2n+7$
except $\lk(yy\cdots yxy,(y\cdots yxy,xy))=\pm 1$, where $\cdots$
represents $(n-2)$~$y$'s, by modifying our example in
Figure~\ref{fig:self-mutation-ex} in a similar way to the construction
of Milnor's link~$L_q$.  The mutant has vanishing $\bar\mu$-invariants
of weight $<2n+7$.  We conjecture that the nontriviality of the
specific $(2n+7)$-linking implies the nontriviality of some
$\bar\mu$-invariant of weight~$(2n+7)$.  The essential difficulty in
the general case is that the technique of the formal Massey products,
which was used in our proof for the simplest case, is not suitable for
large~$n$ since it has exponentially growing computational complexity
in~$n$.  Another interesting question is whether $9$ is the minimal
weight of the $\bar\mu$-invariants not preserved by self-mutation.

\section{Mutation on 3-manifolds}
\label{sec:mutation-3-mfd}

In this section we relate the above results with the effect of
mutation on 3-manifolds.  Let $L$ be a link with $m$ components
in~$S^3$.  By removing an open tubular neighborhood of $L$ and filling
with $m$ solid tori so that each 0-linking longitude bounds a disk in
a solid torus, we obtain a 3-manifold~$M_L$.  We call it the
\emph{surgery manifold} of~$L$.  As expected, the lower central
quotients of surgery manifolds are closely related to the
$\bar\mu$-invariants.  Let $F$ be the free group on $x_1,\ldots,x_m$
as before.

\begin{lem}
\label{lem:surgery-mfd-lcq}
$\pi_1(M_L)/\pi_1(M_L)_q$ is isomorphic to $F/F_q$ if and only if all
$\bar\mu$-invariants of weight~$\le q$ vanish for~$L$.
\end{lem}

\begin{proof}
By Milnor, $G_L/(G_L)_q$ is isomorphic to $F/\langle F_q,
[x_1,w_1],\ldots,[x_m,w_m] \rangle$, where $w_i$ is a word in
$x_1,\ldots,x_m$ representing the $i$-th longitude.  Since
$\pi_1(M_L)$ is the quotient group of $G_L$ modulo the normal subgroup
generated by longitudes, $\pi_1(M_L)/\pi_1(M_L)_q$ is presented as
$F/\langle F_q, w_1,\ldots,w_m \rangle$.

It is well known that all $\bar\mu$-invariants of weight~$\le q$
vanish if and only if $w_1,\ldots,w_m$ are contained in~$F_q$.  If so,
the relations $w_1,\ldots,w_m$ in the above presentation are redundant,
and $\pi_1(M_L)/\pi_1(M_L)_q$ is isomorphic to~$F/F_q$.  Conversely,
if $w_i$ is not contained in $F_q$ for some~$i$, the normal subgroup
$N$ of $F/F_q$ generated by cosets $w_1 F_q,\ldots,w_m F_q$ is
nontrivial.  Since $F/F_q$ is Hopfian (e.g.\
see~\cite{Magnus-Karrass-Solitar:1966-1}),
$\pi_1(M_L)/\pi_1(M_L)_q=(F/F_q)/N$ is not isomorphic to~$F/F_q$.
\end{proof}

The following is an immediate consequence of
Lemma~\ref{lem:surgery-mfd-lcq}.

\begin{thm}
\label{thm:mutative-surgery-mfd-lcq}
Suppose that $L$ is a link with vanishing $\bar\mu$-invariants of
weight~$\le q$, $L^*$ is either a bi-mutant of type $R$ or a
self-mutant of~$L$, and $L^*$ has nonvanishing $\bar\mu$-invariants of
weight~$\le q$.  Then $M_{L^*}$ is obtained by performing mutation
on~$M_L$, but $\pi_1(M_{L^*})/\pi_1(M_{L^*})_q$ is not isomorphic to
$\pi_1(M_L)/\pi_1(M_L)_q$.
\end{thm}

In Sections~\ref{sec:bi-mutation} and~\ref{sec:self-mutation}, we
found examples of bi-mutation of type $R$ and self-mutation satisfying
the hypothesis of Theorem~\ref{thm:mutative-surgery-mfd-lcq}.  Thus
Theorem~\ref{thm:mutation-not-preserving-lcq} follows.

\bibliographystyle{amsplainabbrv}
\bibliography{research}

\providecommand{\bysame}{\leavevmode\hbox to3em{\hrulefill}\thinspace}
\begin{thebibliography}{10}

\bibitem{Cappell-Shaneson:1980-1}
S.~E. Cappell and J.~L. Shaneson, \emph{Link cobordism}, Comment. Math. Helv.
  \textbf{55} (1980), no.~1, 20--49.

\bibitem{Casson-Gordon:1978-1}
A.~Casson and C.~Gordon, \emph{On slice knots in dimension three}, Algebraic
  and geometric topology (Proc. Sympos. Pure Math., Stanford Univ., Stanford,
  Calif., 1976), Part 2, Amer. Math. Soc., Providence, R.I., 1978, pp.~39--53.

\bibitem{Casson-Gordon:1986-1}
\bysame, \emph{Cobordism of classical knots}, \`A la recherche de la topologie
  perdue, Birkh\"auser Boston, Boston, MA, 1986, With an appendix by P. M.
  Gilmer, pp.~181--199.

\bibitem{Cha-Ko:2000-2}
J.~C. Cha and K.~H. Ko, \emph{Signatures of covering links},
  arXiv:math.GT/0108206.

\bibitem{Cha-Ko:2000-1}
\bysame, \emph{Signatures of links in rational homology homology spheres},
  arXiv:math.GT/0108197, to appear in Topology.

\bibitem{Cochran:1985-2}
T.~D. Cochran, \emph{Geometric invariants of link cobordism}, Comment. Math.
  Helv. \textbf{60} (1985), no.~2, 291--311.

\bibitem{Cochran:1987-1}
\bysame, \emph{Link concordance invariants and homotopy theory}, Invent. Math.
  \textbf{90} (1987), no.~3, 635--645.

\bibitem{Cochran:1990-1}
\bysame, \emph{Derivatives of links: {M}ilnor's concordance invariants and
  {M}assey's products}, Mem. Amer. Math. Soc. \textbf{84} (1990), no.~427,
  x+73.

\bibitem{Cochran-Levine:1991-1}
T.~D. Cochran and J.~P. Levine, \emph{Homology boundary links and the
  {A}ndrews-{C}urtis conjecture}, Topology \textbf{30} (1991), no.~2, 231--239.

\bibitem{Cochran-Orr:1993-1}
T.~D. Cochran and K.~E. Orr, \emph{Not all links are concordant to boundary
  links}, Ann. of Math. (2) \textbf{138} (1993), no.~3, 519--554.

\bibitem{Cochran-Orr:1994-1}
\bysame, \emph{Homology boundary links and {B}lanchfield forms: concordance
  classification and new tangle-theoretic constructions}, Topology \textbf{33}
  (1994), no.~3, 397--427.

\bibitem{Cochran-Orr-Teichner:1999-1}
T.~D. Cochran, K.~E. Orr, and P.~Teichner, \emph{Knot concordance, whitney
  towers and ${L}^2$-signatures}, arXiv:math.GT/9908117, 1999.

\bibitem{Duval:1986-1}
J.~Duval, \emph{Forme de {B}lanchfield et cobordisme d'entrelacs bords},
  Comment. Math. Helv. \textbf{61} (1986), no.~4, 617--635.

\bibitem{Kaiser:1992-1}
U.~Kaiser, \emph{Homology boundary links and fusion constructions}, Osaka J.
  Math. \textbf{29} (1992), no.~3, 573--593.

\bibitem{Kania-Bartoszynska:1993-1}
J.~Kania-Bartoszynska, \emph{Examples of different $3$-manifolds with the same
  invariants of {W}itten and {R}eshetikhin-{T}uraev}, Topology \textbf{32}
  (1993), no.~1, 47--54.

\bibitem{Kawauchi:1994-1}
A.~Kawauchi, \emph{Topological imitation, mutation and the quantum {S}{U}$(2)$
  invariants}, J. Knot Theory Ramifications \textbf{3} (1994), no.~1, 25--39.

\bibitem{Kawauchi:1996-2}
\bysame, \emph{Mutative hyperbolic homology $3$-spheres with the same {F}loer
  homology}, Geom. Dedicata \textbf{61} (1996), no.~2, 205--217.

\bibitem{Kirk:1989-1}
P.~Kirk, \emph{Mutations of homology spheres and {C}asson's invariant}, Math.
  Proc. Cambridge Philos. Soc. \textbf{105} (1989), no.~2, 313--318.

\bibitem{Kirk-Livingston:1999-1}
P.~Kirk and C.~Livingston, \emph{Concordance and mutation},
  arXiv:math.GT/9912174.

\bibitem{Kirk-Livingston:1999-3}
P.~Kirk and C.~Livingston, \emph{Twisted knot polynomials: inversion, mutation
  and concordance}, Topology \textbf{38} (1999), no.~3, 663--671.

\bibitem{Ko:1987-1}
K.~H. Ko, \emph{Seifert matrices and boundary link cobordisms}, Trans. Amer.
  Math. Soc. \textbf{299} (1987), no.~2, 657--681.

\bibitem{Krushkal:1998-1}
V.~S. Krushkal, \emph{Additivity properties of {M}ilnor's
  $\overline\mu$-invariants}, J. Knot Theory Ramifications \textbf{7} (1998),
  no.~5, 625--637.

\bibitem{Levine:1989-1}
J.~P. Levine, \emph{Link concordance and algebraic closure. {I}{I}}, Invent.
  Math. \textbf{96} (1989), no.~3, 571--592.

\bibitem{Levine:1989-2}
\bysame, \emph{Link concordance and algebraic closure of groups}, Comment.
  Math. Helv. \textbf{64} (1989), no.~2, 236--255.

\bibitem{Lickorish:1993-1}
W.~B.~R. Lickorish, \emph{Distinct $3$-manifolds with all {S}{U}$(2)\sb q$
  invariants the same}, Proc. Amer. Math. Soc. \textbf{117} (1993), no.~1,
  285--292.

\bibitem{Magnus-Karrass-Solitar:1966-1}
W.~Magnus, A.~Karrass, and D.~Solitar, \emph{Combinatorial group theory:
  {P}resentations of groups in terms of generators and relations}, Interscience
  Publishers [John Wiley \& Sons, Inc.], New York-London-Sydney, 1966.

\bibitem{Meyerhoff-Ruberman:1990-1}
R.~Meyerhoff and D.~Ruberman, \emph{Mutation and the $\eta$-invariant}, J.
  Differential Geom. \textbf{31} (1990), no.~1, 101--130.

\bibitem{Milnor:1954-1}
J.~W. Milnor, \emph{Link groups}, Ann. of Math. (2) \textbf{59} (1954),
  177--195.

\bibitem{Milnor:1957-1}
\bysame, \emph{Isotopy of links. {A}lgebraic geometry and topology}, A
  symposium in honor of S. Lefschetz, Princeton University Press, Princeton, N.
  J., 1957, pp.~280--306.

\bibitem{Mio:1987-1}
W.~Mio, \emph{On boundary-link cobordism}, Math. Proc. Cambridge Philos. Soc.
  \textbf{101} (1987), no.~2, 259--266.

\bibitem{Morton-Traczyk:1988-1}
H.~R. Morton and P.~Traczyk, \emph{The {J}ones polynomial of satellite links
  around mutants}, Braids (Santa Cruz, CA, 1986), Amer. Math. Soc., Providence,
  RI, 1988, pp.~587--592.

\bibitem{Orr:1989-1}
K.~E. Orr, \emph{Homotopy invariants of links}, Invent. Math. \textbf{95}
  (1989), no.~2, 379--394.

\bibitem{Porter:1980-1}
R.~Porter, \emph{Milnor's $\bar \mu $-invariants and {M}assey products}, Trans.
  Amer. Math. Soc. \textbf{257} (1980), no.~1, 39--71.

\bibitem{Porti:1997-1}
J.~Porti, \emph{Torsion de {R}eidemeister pour les vari\'et\'es hyperboliques},
  Mem. Amer. Math. Soc. \textbf{128} (1997), no.~612, x+139.

\bibitem{Rong:1994-1}
Y.~W. Rong, \emph{Mutation and {W}itten invariants}, Topology \textbf{33}
  (1994), no.~3, 499--507.

\bibitem{Ruberman:1987-1}
D.~Ruberman, \emph{Mutation and volumes of knots in ${S}\sp 3$}, Invent. Math.
  \textbf{90} (1987), no.~1, 189--215.

\bibitem{Ruberman-1999:1}
\bysame, \emph{Mutation and gauge theory. {I}. {Y}ang-{M}ills invariants},
  Comment. Math. Helv. \textbf{74} (1999), no.~4, 615--641.

\bibitem{Sato:1984-1}
N.~Sato, \emph{Cobordisms of semiboundary links}, Topology Appl. \textbf{18}
  (1984), no.~2-3, 225--234.

\bibitem{Stallings:1965-1}
J.~Stallings, \emph{Homology and central series of groups}, J. Algebra
  \textbf{2} (1965), 170--181.

\bibitem{Stein:1990-1}
D.~Stein, \emph{Massey products in the cohomology of groups with applications
  to link theory}, Trans. Amer. Math. Soc. \textbf{318} (1990), no.~1,
  301--325.

\bibitem{Turaev:1976-1}
V.~G. Turaev, \emph{The {M}ilnor invariants and {M}assey products}, Zap. Nau\v
  cn. Sem. Leningrad. Otdel. Mat. Inst. Steklov. (LOMI) \textbf{66} (1976),
  189--203, 209--210, Studies in topology, II.

\end{thebibliography}

\end{document}